\documentclass[12pt]{article}

\usepackage[T1]{fontenc}
\usepackage[utf8]{inputenc}
\usepackage[english]{babel}
\usepackage{amsmath,amsfonts,amsthm,amssymb,MnSymbol}
\usepackage{tikz}
\usepackage{tikz-qtree}
\usepackage{hyperref,url}
\usepackage{bussproofs}
\usepackage{xcolor}
\usepackage[margin=2.5cm]{geometry}
\usepackage{enumitem}
\hypersetup{colorlinks=true,linkcolor=blue,citecolor=blue,urlcolor=blue}
\usepackage{cite}

\DeclareMathOperator{\cd}{\mathit{cd}}

\EnableBpAbbreviations

\newtheorem{thm}{Theorem}
\newtheorem{lem}{Lemma}
\newtheorem{prop}{Proposition}

\theoremstyle{definition}

\newcommand{\Rto}{\Rightarrow}

\providecommand{\G}{\Gamma}
\renewcommand{\G}{\Gamma}
\newcommand{\D}{\Delta}

\newcommand{\qgl}{\mathsf{QGL}}
\newcommand{\QGL}{\mathsf{QGL}}

\title{Topological Completeness of the Predicate Modal Calculus $\QGL$ with Non-Well-Founded Derivations}

\author{Pavel Razumnyi\\
\normalsize{\textit{HSE University}}\\
\normalsize{\textit{Faculty of Mathematics}}\\
\normalsize{\textit{6 Usachyova St, Moscow, 119048, Russia}}\\
\normalsize{\textit{razumnyyp@gmail.com}}\\\\ Daniyar Shamkanov\\
\normalsize{\textit{Steklov Mathematical Institute}}\\
\normalsize{\textit{of the Russian Academy of Sciences}}\\
\normalsize{\textit{8 Gubkina St, Moscow, 119333, Russia}}\\
\normalsize{\textit{daniyar.shamkanov@gmail.com}}}

\date{}

\begin{document}
\maketitle

\sloppy

\begin{abstract}
According to the result of F. Montagna \cite{montagna1984}, the predicate version of the G\"odel--L\"ob modal logic $\mathsf{GL}$, denoted by $\mathsf{QGL}$, is incomplete with respect to relational semantics. In this paper, we study an extension of $\mathsf{QGL}$ by means of non-well-founded derivations and establish completeness of the resulting infinitary system with respect to topological semantics.
\end{abstract}

\section{Introduction}
\label{s0}
The G\"odel--L\"ob provability logic $\mathsf{GL}$ is a propositional modal system widely known for being sound and complete with respect to the arithmetical interpretation, under which the modal connective $\Box$ is understood as the standard provability predicate $\mathsf{Prov}_{\mathsf{PA}}(\cdot)$ of Peano arithmetic $\mathsf{PA}$ \cite{solovay1976}. In other words, $\mathsf{GL}$ describes the properties of the provability predicate expressible in the propositional modal language and establishable by means of Peano arithmetic.
$\mathsf{GL}$ is also complete with respect to its relational semantics \cite{segerberg1971}; however, completeness holds only in the weak sense.
With respect to the more general topological, or neighborhood, semantics, $\mathsf{GL}$ is strongly complete \cite{shehtman2005}.

It was natural to assume that the predicate provability logic of Peano arithmetic is given by the calculus $\mathsf{QGL}$, the predicate version of $\mathsf{GL}$. However, F. Montagna \cite{montagna1984} established that the calculus $\mathsf{QGL}$ is complete neither with respect to relational semantics nor with respect to the arithmetical one. Moreover, as shown by V. A. Vardanyan \cite{vardanyan1986}, the predicate provability logic of Peano arithmetic is $\Pi^0_2$-complete (hence, this logic has no effective axiomatization). Let us also recall the theorem of S. N. Artemov \cite{artemov1985}, according to which the set of predicate modal formulas true in the standard model under any arithmetical interpretation is non-arithmetical.

Several years ago, D. S. Shamkanov considered cyclic and non-well-founded derivations in the logic $\mathsf{GL}$ \cite{shamkanov2014, shamkanov2020}. It is worth recalling that a \emph{cyclic derivation} can be obtained from an ordinary formal derivation from hypotheses (a finite tree of formulas built according to the inference rules of a given calculus)
by successively applying the following operation: we identify a certain leaf with an inner vertex of the derivation lying on the path between the root and this leaf and labeled with the same formula as the leaf, by means of a so-called back reference. In this case, special conditions are usually imposed on these references to ensure the correctness of the derivation.
\emph{Non-well-founded derivations} are a generalization of cyclic ones and represent infinite trees of formulas built according to the inference rules, on whose infinite branches special conditions are also imposed to ensure correctness. In the case of $\mathsf{GL}$, every infinite branch must contain infinitely many applications of the rule $(\mathsf{nec})$, i.e., the necessitation rule. In \cite{shamkanov2020}, D. S. Shamkanov considered the logic $\mathsf{GL}$ extended by non-well-founded derivations and proved its completeness with respect to algebraic and topological semantics.

In the present paper, we study
the extension of the predicate calculus $\mathsf{QGL}$ by means of non-well-founded derivations, which we denote by $\mathsf{QGL}_\infty$. The main contribution of the paper is the proof of soundness and completeness results for the calculus $\mathsf{QGL}_\infty$ with respect to topological semantics (given by the class of predicate topological $\mathsf{GL}$-frames with constant domains). The completeness proof relies on an auxiliary sequent calculus admitting non-well-founded proofs and employs K. Sch\"utte's method of reduction trees.

The rest of the paper is organized as follows. In \S \ref{s1}, we introduce non-well-founded derivations in the calculus $\mathsf{QGL}$ and define the calculus $\mathsf{QGL}_\infty$. In \S \ref{s2}, we define topological semantics and establish the topological soundness of the system $\mathsf{QGL}_\infty$. In the next section, we introduce a non-well-founded sequent calculus and prove its topological completeness.
In \S \ref{s4}, we establish the topological completeness of the calculus $\mathsf{QGL}_\infty$. In the concluding section, we outline directions for future research.

\section{Non-well-founded derivations}
\label{s1}

In this section, we recall how the predicate version of the G\"odel--L\"ob modal logic $\mathsf{GL}$, denoted by $\mathsf{QGL}$, is defined, and we introduce non-well-founded derivations in $\mathsf{QGL}$.

To simplify further considerations, we work with a signature containing no function symbols or constants and build $\mathsf{QGL}$-formulas from atomic formulas $P(y_0,\dotsc, y_{n-1})$ and the formula $\bot$ using the propositional connective $\to$, the modal connective $\Box$, and the universal quantifier $\forall$. The remaining Boolean connectives, the modal connective $\Diamond$, and the existential quantifier $\exists$ are understood as abbreviations:
\begin{gather*}
\neg A \coloneq A\to \bot,\qquad\top \coloneq \neg \bot,\qquad A\wedge B \coloneq \neg (A\to \neg B), \qquad A\vee B \coloneq \neg A\to B,\\  A \leftrightarrow B\coloneq(A\to B)\wedge (B \to A),\qquad\Diamond A \coloneq \neg\Box \neg A, \qquad \exists x\: A\coloneq \neg \forall x\: \neg A.
\end{gather*}
The set of individual variables of the given language $\{x_0, x_1, x_2, \dotsc \}$ is denoted by the symbol $\mathit{Var}$, and the set of variables having a free occurrence in a formula $A$ by $\mathit{FV}(A)$. Here, $\mathit{Fm}_\Omega$ is the set of all formulas of a fixed signature $\Omega$.

In what follows, we denote by $A[x\mapsto y]$ the formula obtained from a formula $A$ by replacing all free occurrences of the variable $x$ with $y$ (renaming bound occurrences of variables if necessary). The formula $A[x\mapsto y]$ is defined recursively by the following conditions:
\begin{gather*}\bot [x\mapsto y] =\bot, \qquad \qquad \quad (\Box B)[x\mapsto y]= \Box B[x\mapsto y],\\
(B\to C)[x\mapsto y] = B [x\mapsto y] \to C [x\mapsto y].
\end{gather*}
Moreover,
$P(z_0,\dotsc, z_{n-1})[x\mapsto y]= P(z_0[x\mapsto y],\dotsc, z_{n-1}[x\mapsto y])$, where
\begin{gather*}
z[x\mapsto y]=
\begin{cases}
y , & \text{if $z = x$},\\
z, & \text{otherwise}.
\end{cases}
\end{gather*}
Also,
\begin{gather*}
(\forall z\: B)[x\mapsto y]=
\begin{cases}
\forall z\: B , & \text{if $z=x$},\\
\forall z\: B[x\mapsto y], & \text{if $z\neq x$ and $z\neq y$},\\
\forall x_j\: (B[z\mapsto x_j])[x\mapsto y], & \text{if $z\neq x$ and $z= y$},
\end{cases}
\end{gather*}
where $j = \min
\{i\in \mathbb{N}\mid
\text{$x_i\notin \mathit{FV}(B)$, $x_i\neq x$ and $x_i \neq y$}\}$.



The following axioms and inference rules define the calculus $\mathsf{QGL}$.\medskip

\textit{Axioms:}
\begin{itemize}
\item[(i)] tautologies of classical propositional logic,
\item[(ii)] $\Box (A\to B)\to (\Box A \to \Box B)$,
\item[(iii)] $\Box (\Box A \to A) \to  \Box A$,
\item[(iv)] $\forall x\: A \to A[x\mapsto y]$,
\item[(v)] $\forall x\: (A\to B) \to (A\to \forall x\: B)$, where $x\notin \mathit{FV}(A)$.
\end{itemize}

\textit{Inference rules:}
\[
\AXC{$A$}
\AXC{$ A\to B$}
\LeftLabel{$\mathsf{mp}$}
\RightLabel{, }
\BIC{$ B$}
\DisplayProof
\qquad
\AXC{$A$}
\LeftLabel{$\mathsf{nec}$}
\RightLabel{,}
\UIC{$\Box A$}
\DisplayProof
\qquad
\AXC{$A$}
\LeftLabel{$\mathsf{gen}$}
\RightLabel{.}
\UIC{$\forall x\: A$}
\DisplayProof
\]

The first three items of the list of axioms of the system $\mathsf{QGL}$ repeat the corresponding items of the axiomatic calculus for the logic $\mathsf{GL}$. The last two items are the axioms for the universal quantifier of classical predicate calculus in a language without function symbols or constants. Note that the transitivity principle $\Box A \to \Box \Box A$ is provable in the calculus $\mathsf{QGL}$ exactly as in $\mathsf{GL}$.

Now we define infinitary derivations in the calculus $\mathsf{QGL}$. By a \emph{non-well-founded derivation}, or \emph{$\infty$-derivation}, we mean a (possibly infinite) tree whose vertices are labeled with $\mathsf{QGL}$-formulas, and the transition from children to parents is performed according to the inference rules ($\mathsf{mp}$), ($\mathsf{nec}$), and ($\mathsf{gen}$). In addition, we require that every infinite branch of an $\infty$-derivation intersects the rule ($\mathsf{nec}$) infinitely many times. A non-well-founded derivation all of whose leaves are labeled with axioms of $\mathsf{QGL}$ will be called a \emph{non-well-founded proof}, or \emph{$\infty$-proof}.


We give an example of an $\infty$-derivation, assuming that for each $n\in \mathbb{N}$ the signature $\Omega$ contains a unary predicate symbol $P_n$:
\begin{gather*}\label{example1}
\AXC{$\vdots$}
\noLine
\UIC{$\Box  \forall x_2\: P_2(x_2)$}
\AXC{$\Box \forall x_2\: P_2(x_2) \to P_1(x_1)$}
\LeftLabel{$\mathsf{mp}$}
\BIC{$  P_1(x_1)$}
\LeftLabel{$\mathsf{gen}$}
\UIC{$\forall x_1\: P_1(x_1)$}
\LeftLabel{$\mathsf{nec}$}
\UIC{$\Box \forall x_1\: P_1(x_1)$}
\AXC{$\Box \forall x_1\: P_1(x_1) \to P_0(x_0)$}
\LeftLabel{$\mathsf{mp}$}
\RightLabel{ .}
\BIC{$P_0(x_0)$}
\DisplayProof
\end{gather*}
In this tree there is a unique infinite branch, which intersects the rule ($\mathsf{nec}$) infinitely many times, as required.

\emph{The main fragment of an $\infty$-derivation $\delta$} is the tree obtained from $\delta$ by cutting each branch at the application of the rule ($\mathsf{nec}$) closest to the root. \emph{The local height of $\delta$,} denoted by $\lvert \delta \rvert$, is defined as the length of the longest branch of the main fragment of $\delta$. If $\delta$ consists of a single vertex, then $\lvert \delta \rvert=0$.

Note that the local height of the above example is equal to $1$. Moreover, its main fragment has the form:
\begin{gather*}
\AXC{$\Box\forall x_1\: P_1(x_1)\qquad\qquad\quad \;\;\;$}
\AXC{$\Box \forall x_1\: P_1(x_1) \to P_0(x_0)$}
\LeftLabel{$\mathsf{mp}$}
\RightLabel{ .}
\BIC{$P_0(x_0)$}
\DisplayProof
\end{gather*}


We denote by $\mathsf{QGL}_\infty$ the extension of $\mathsf{QGL}$ in which non-well-founded derivations are allowed. Thus, $\mathsf{QGL}_\infty\vdash A$ means that there exists an $\infty$-proof whose root is labeled with the formula $A$.


\section{Topological semantics}
\label{s2}

In the first half of this section, we recall the main properties of scattered topological spaces. Note that the class of these spaces corresponds exactly to the G\"odel--L\"ob modal logic $\mathsf{GL}$ under the topological interpretation (for a detailed overview of the topological semantics of $\mathsf{GL}$, we refer the reader to \cite{beklemishev2014}). In the second half of the section, we define the topological semantics of the system $\mathsf{QGL}$ and show that, with respect to this interpretation, all theorems of $\mathsf{QGL}_\infty$ are valid.


Recall that a point $w$ in a topological space is \emph{isolated} if $\{w\}$ is open. A topological space $(X, \tau)$ is called \emph{scattered} if every nonempty subset of it (as a space with the induced topology) contains an isolated point. Every ordinal with the interval topology gives an example of such a space. The ordinal $0$ also gives us such an example, since we allow a topological space to be empty, i.e., equal to $(\emptyset, \{\emptyset\})$.

\emph{A neighborhood of a point $w$} in a topological space $(X, \tau)$ is any open set $U$ containing $w$. A set $U$ is called a \emph{punctured neighborhood of $w$} if $w \notin U$ and $U\cup \{w\}$ is open.
For $Y\subset X$, the \emph{derived set $d_\tau (Y)$} is defined as the set of limit points of $Y$:
\[w \in d_\tau (Y) \Longleftrightarrow \forall U \in \tau \; (w\in U \Rightarrow \exists w^\prime \neq w \; (w^\prime \in U \cap Y))\;.\]
The \emph{co-derived set $\cd_\tau(Y)$} is $X \setminus d_\tau (X \setminus Y)$.
From this definition it follows that $w \in\cd_\tau(Y)$ if and only if some punctured neighborhood of $w$ is contained in $Y$. It is easy to see that $\cd_\tau$ is a monotone operation on $\mathcal{P}(X)$, i.e., $\cd_\tau(Y_0)\subset \cd_\tau(Y_1)$ whenever $Y_0\subset Y_1\subset X$. Note also that a set $Y$ is open if and only if $Y\subset \cd_\tau (Y)$.

For an arbitrary topological space $(X, \tau)$, consider the transfinite sequence $\cd^{\alpha}_\tau (\emptyset)$ given by the recursive condition
\[\cd^{\alpha}_\tau (\emptyset) = \bigcup_{\beta < \alpha} \cd_\tau(\cd^{\beta}_\tau (\emptyset)).\]
\begin{lem}
In any topological space $(X, \tau)$, the following hold:
\begin{itemize}
\item $\cd^{\beta}_\tau (\emptyset) \subset \cd^{\alpha}_\tau (\emptyset)$, if $\beta\leqslant\alpha$;
\item $\cd^{\alpha+1}_\tau (\emptyset)= \cd_\tau(\cd^{\alpha}_\tau (\emptyset))$;
\item $\cd^{\alpha}_\tau (\emptyset)$ is open for every $\alpha$.
\end{itemize}
\end{lem}
\begin{proof}
It follows directly from the definition that the sequence under consideration is non-strictly increasing, i.e., $\cd^{\beta}_\tau (\emptyset) \subset \cd^{\alpha}_\tau (\emptyset)$ whenever $\beta\leqslant\alpha$.

Moreover, if $\beta\leqslant\alpha$, then $ \cd_\tau(\cd^{\beta}_\tau (\emptyset))\subset \cd_\tau(\cd^{\alpha}_\tau (\emptyset))$ by the monotonicity of the operation $\cd_\tau$. Hence we obtain
\[\cd_\tau(\cd^{\alpha}_\tau (\emptyset)) \subset \bigcup_{\beta < \alpha+1} \cd_\tau(\cd^{\beta}_\tau (\emptyset))\subset \cd_\tau(\cd^{\alpha}_\tau (\emptyset))\]
for every ordinal $\alpha$. Consequently, $\cd^{\alpha+1}_\tau (\emptyset)= \cd_\tau(\cd^{\alpha}_\tau (\emptyset))$.

Finally, $\cd^{\alpha}_\tau (\emptyset) \subset \cd^{\alpha+1}_\tau (\emptyset)= \cd_\tau(\cd^{\alpha}_\tau (\emptyset))$, i.e., $\cd^{\alpha}_\tau (\emptyset)$ is contained in its co-derived set. Consequently, $\cd^{\alpha}_\tau (\emptyset)$ is open for every $\alpha$.
\end{proof}

\begin{prop}[G. Cantor]\label{Cantor}
A topological space $(X, \tau)$ is scattered if and only if $\cd^\alpha_\tau (\emptyset) = X$ for some $\alpha$.
\end{prop}
For a scattered space $(X, \tau)$, the \emph{rank $\rho_\tau(w)$ of a point $w$} is defined as the least ordinal $\beta$ such that $w \in cd^{\beta+1}_\tau(\emptyset)$. It follows from the previous proposition that the rank function is defined at every point. Moreover, since
\[\cd^{\alpha}_\tau (\emptyset)=  \bigcup_{\beta < \alpha} \cd^{\beta+1}_\tau (\emptyset),\]
it is obvious that $\cd^{\alpha}_\tau (\emptyset) =  \{w\in X \mid \rho_\tau (w)<\alpha\}$.

\begin{lem}\label{punctured neighbourhood}
For every point $w$ of a scattered space $(X, \tau)$, the set $\{w^\prime\in X \mid \rho_\tau(w^\prime) < \rho_\tau(w)\}$ is a punctured neighborhood of $w$.
\end{lem}
\begin{proof}
Let us check that $\{w\} \cup\{w^\prime\in X \mid \rho_\tau(w^\prime) < \rho_\tau(w)\}$ is contained in its co-derived set. Note that $\{w\} \subset  \cd^{\rho_\tau (w) + 1}_\tau (\emptyset)$ and
\[
\{w^\prime\in X \mid \rho_\tau(w^\prime) < \rho_\tau(w)\}= \cd^{\rho_\tau (w)}_\tau (\emptyset) \subset cd^{\rho_\tau (w) + 1}_\tau (\emptyset) .
\]
Thus,
\begin{multline*}
\{w\} \cup\{w^\prime\in X \mid \rho_\tau(w^\prime) < \rho_\tau(w)\} \subset  \cd^{\rho_\tau (w) + 1}_\tau (\emptyset) = \cd_\tau (\cd^{\rho_\tau (w) }_\tau (\emptyset)) \subset\\
\subset  \cd_\tau (\{w\} \cup \cd^{\rho_\tau (w) }_\tau (\emptyset))   = cd_\tau (\{w\} \cup\{w^\prime\in X \mid \rho_\tau(w^\prime) < \rho_\tau(w)\} ).
\end{multline*}
Consequently, the set $\{w\} \cup\{w^\prime\in X \mid \rho_\tau(w^\prime) < \rho_\tau(w)\}$ is open. Since $w\notin \{w^\prime\in X \mid \rho_\tau(w^\prime) < \rho_\tau(w)\}$, we obtain that $\{w^\prime\in X \mid \rho_\tau(w^\prime) < \rho_\tau(w)\}$ is a punctured neighborhood of $w$.

\end{proof}

Now we proceed to the definition of the topological semantics of the system $\mathsf{QGL}$. By a \emph{predicate topological $\mathsf{GL}$-frame} we mean a tuple $(X,\tau, D)$, where $(X,\tau)$ is a scattered topological space and $D$ is a nonempty set called the \emph{domain} of the frame. We also call predicate topological $\mathsf{GL}$-frames \emph{$\mathsf{QGL}$-frames}.

A \emph{valuation, or interpretation, in $D$} is a function assigning to each predicate symbol of the signature of arity $n$ a relation on the set $D$ of the same arity. We call a \emph{predicate topological $\mathsf{GL}$-model}, or \emph{$\mathsf{QGL}$-model}, any tuple $\mathcal{M}=(X,\tau, D, \xi)$ consisting of a $\mathsf{QGL}$-frame $(X,\tau,D)$ and an indexed family $\xi=(\xi_w)_{w\in X}$ of valuations in $D$. The elements of the set $X$ are usually called \emph{worlds} of the model (or of the corresponding frame).

A \emph{variable assignment in $D$} is a function from the set of individual variables $\mathit{Var}=\{x_0, x_1, x_2,\dotsc\}$ to the domain $D$. The \emph{truth of a formula $A$ at a world $w$ of a model $\mathcal{M}$ under an assignment $h$} is defined as follows:
\begin{itemize}
\item $\mathcal{M},w,h \nvDash\bot$,
\item $\mathcal{M},w,h \vDash P(y_0,\dotsc, y_{n-1}) \Longleftrightarrow (h(y_0),\dotsc,h(y_{n-1}))\in \xi_w(P)$,
\item $\mathcal{M},w,h \vDash A\to B \Longleftrightarrow \mathcal{M},w,h \nvDash A \text{ or } \mathcal{M},w,h \vDash B$,
\item $\mathcal{M},w,h \vDash \Box A \Longleftrightarrow \exists U\in \tau\: (w\in U\text{ and } \forall w^\prime\in U\setminus \{w\}\;\mathcal{M},w^\prime,h \vDash A)$,
\item $\mathcal{M},w,h \vDash \forall x\: A \Longleftrightarrow \forall a\in D\: \mathcal{M},w,h^a_x \vDash A$,
\end{itemize}
where $h^a_x$ is the assignment such that $h^a_x (x)=a$ and $ h^a_x (y)=h(y)$ for all $y\in \mathit{Var}\setminus \{x\}$. A \emph{formula $A$ is valid} if for every $\mathsf{QGL}$-model $\mathcal{M}$ the formula $A$ is true at every world of $\mathcal{M}$ under every variable assignment.


It is worth noting that the Barcan formula $\forall y \Box P(y) \to \Box \forall y\: P(y)$, being valid in relational semantics with constant domains, is not valid in the topological semantics of the system $\mathsf{QGL}$. To show this, suppose that the signature $\Omega$ contains a unary predicate symbol $P$, and consider the $\mathsf{QGL}$-model $\mathcal{M}^\prime$
\begin{center}
\begin{tikzpicture}[grow'=up]
\Tree [.$r$ $w_0$ $w_1$ $w_2$ $\dots$ ]
\end{tikzpicture}
\end{center}
with the set of worlds $X^\prime= \{ r \}\cup \{ w_0, w_1, w_2, \dotsc \}$, the domain $D^\prime=\{a_0,a_1,a_2, \dotsc\}$, and an arbitrary family of valuations $\xi^\prime$ such that $\xi^\prime_{w_n}(P)= \{a_k \mid k< n\}$.
Here, a set of worlds $U\in \tau^\prime$ if and only if
\[r\in U \Longrightarrow\exists k\in \mathbb{N}\: \forall n> k \: w_n \in U.\]
It is obvious that $(X^\prime, \tau^\prime)$ is a scattered topological space and $\mathcal{M}^\prime=(X^\prime,\tau^\prime, D^\prime, \xi^\prime)$ is indeed a $\mathsf{QGL}$-model. Furthermore, for the assignment $h \colon x_i \mapsto a_i$, we have $\mathcal{M}^\prime, r , h \vDash \forall  y \Box P(y) $, since $\mathcal{M}^\prime, w_n , h^{a_k}_y \vDash P(y) $ whenever $n> k$, and $\mathcal{M}^\prime, r , h^{a_k}_y \vDash \Box P(y) $ for every $k\in \mathbb{N}$. At the same time, $\mathcal{M}^\prime, r , h \nvDash  \Box \forall  y\: P(y) $, since $\mathcal{M}^\prime, w_n , h^{a_{n}}_y \nvDash P(y) $ and $\mathcal{M}^\prime, w_n , h \nvDash   \forall  y\: P(y) $ for each $n\in \mathbb{N}$. Thus, $\mathcal{M}^\prime, r , h \nvDash \forall  y \Box P(y) \to \Box \forall  y\: P(y) $.

Let us analyze one more example. Although the formula $\forall x\:\exists y\: \Box (\Box P(y) \to P(x))\to \forall x\: \Box P(x)$ is valid in the class of predicate relational $\mathsf{GL}$-frames with expanding domains (see \cite{montagna1984}), it is not valid in the topological semantics of the system $\mathsf{QGL}$. Indeed, suppose that the signature $\Omega$ contains a unary predicate symbol $P$, and consider the $\mathsf{QGL}$-model $\mathcal{M}^{\prime\prime}$
\begin{center}
\begin{tikzpicture}[grow'=up]
\Tree [.$r$ [.$w_{0,0}$ ] [.$w_{1,1}$ $w_{1,0}$ ] [.$w_{2,2}$ [.$w_{2,1}$ $w_{2,0}$ ] ] [.$\dots$ ]  ]
\end{tikzpicture}
\end{center}
with the set of worlds $X^{\prime\prime}= \{ r \}\cup \{w_{n,m}\mid n\in \mathbb{N}, m\in \mathbb{N} \text{ and } m\leqslant n\}$, the domain $D^{\prime\prime} =\{a_0,a_1,a_2, \dotsc\}$, and an arbitrary family of valuations $\xi^{\prime\prime}$ such that $\xi^{\prime\prime}_{w_{n,m}}(P)= \{a_k \mid k< n-m\}$.
Here, the topology $ \tau^{\prime\prime}$ on $X^{\prime\prime}$ consists of all sets $U\subset X^{\prime\prime} $ for which
\begin{gather*}
r\in U \Longrightarrow\exists k\in \mathbb{N}\: \forall n> k \: \forall m\leqslant n \: w_{n,m} \in U,\qquad
w_{n,m} \in U \Longrightarrow \forall k\leqslant m \: w_{n,k} \in U.
\end{gather*}
We see that $(X^{\prime\prime}, \tau^{\prime\prime})$ is a scattered topological space and $\mathcal{M}^{\prime\prime}=(X^{\prime\prime},\tau^{\prime\prime}, D^{\prime\prime}, \xi^{\prime\prime})$ is a $\mathsf{QGL}$-model. Also, for the assignment $h \colon x_i \mapsto a_i$, we have $\mathcal{M}^{\prime\prime}, r , h \vDash \forall x\:\exists y\: \Box (\Box P(y) \to P(x)) $ because $\mathcal{M}^{\prime\prime}, w_{n,m} , (h^{a_k}_x)^{a_{k+1}}_y  \vDash \Box P(y) \to P(x) $ whenever $n> k$ and $m\leqslant n$, and $\mathcal{M}^{\prime\prime}, r , (h^{a_k}_x)^{a_{k+1}}_y \vDash \Box (\Box P(y) \to P(x)) $ for every $k\in \mathbb{N}$.
At the same time, $\mathcal{M}^{\prime\prime}, r , h \nvDash  \forall x\: \Box P(x) $, since $\mathcal{M}^{\prime\prime}, w_{n,n} , h^{a_{0}}_x \nvDash P(x) $ for each $n\in \mathbb{N}$ and $\mathcal{M}^{\prime\prime}, r , h^{a_{0}}_x \nvDash \Box P(x)$.
Thus, $\mathcal{M}^{\prime\prime}, r , h \nvDash \forall x\:\exists y\: \Box (\Box P(y) \to P(x))\to \forall x\: \Box P(x) $. We have completed the analysis of the example, showing that the formula $\forall x\:\exists y\: \Box (\Box P(y) \to P(x))\to \forall x\: \Box P(x)$ is not valid in the topological semantics of the system $\mathsf{QGL}$.

The following lemma is established by a straightforward induction on the construction of the formula $A$, so we omit the proof.

\begin{lem}\label{submodel}
Let $\mathcal{M'}=(X',\tau', D, \xi^\prime)$ and $\mathcal{M}=(X,\tau, D, \xi)$ be $\qgl$-models such that $X'$ is an open subset of $X$, $\tau^\prime$ is the topology on $X'$ induced by the topology $\tau$, and $\xi^\prime_w =\xi_w$ for all $w\in X^\prime$. Then for every formula $A$, every world $w \in X^\prime$, and every assignment $h$, the following equivalence holds:
\[\mathcal{M'}, w, h \vDash A \Longleftrightarrow \mathcal{M}, w, h \vDash A .
\]
\end{lem}

We conclude the section with the statement of topological soundness of the system $\mathsf{QGL}_\infty$.
\begin{prop}[topological soundness]\label{soundness}
If $\mathsf{QGL}_\infty \vdash A$, then $A$ is valid.
\end{prop}
\begin{proof}
Consider an arbitrary $\mathsf{QGL}$-model $\mathcal{M}= (X, \tau, D, \xi)$, a world $w$ of this model, and a variable assignment $h\colon \mathit{Var}\to D$. Also, consider an arbitrary $\infty$-proof $\pi$.
Using transfinite induction on $\rho_\tau(w)$ together with a nested induction on $\lvert \pi \rvert$, we prove that $\mathcal{M}, w, h \vDash A$, where $A$ is the formula standing at the root of the $\infty$-proof $\pi$.

Case 1. If the $\infty$-proof $\pi$ consists of a single vertex, then $A$ is an axiom of $\mathsf{QGL}$. Let us verify that $\mathcal{M}, w,h \vDash A$ in the case when $A$ has the form $\Box (\Box B \to B) \to  \Box B$. The cases of all the other axioms are considerably simpler, so we will not consider them.

By the induction hypothesis, we have
\begin{gather}\label{form1}
\forall w^\prime\in X\; (\rho_\tau(w^\prime) < \rho_\tau(w) \Longrightarrow \mathcal{M}, w^\prime, h \vDash \Box (\Box B \to B) \to  \Box B).
\end{gather}
If $\mathcal{M}, w,h \nvDash \Box (\Box B \to B) $, then $\mathcal{M}, w^\prime, h \vDash \Box (\Box B \to B) \to  \Box B$, i.e., $\mathcal{M}, w,h \vDash A $. Now suppose that $\mathcal{M}, w,h \vDash \Box (\Box B \to B) $. By the definition of the truth relation for a $\mathsf{QGL}$-model, there is a neighborhood $U$ of $w$ such that $\mathcal{M}, w^\prime,h \vDash \Box B \to B $ for all $w^\prime\in U\setminus \{w\} $. If $\mathcal{M}, w,h \nvDash \Box B \to B $, then $\mathcal{M}, w,h \vDash \Box B  $ and $\mathcal{M}, w,h \vDash A $. If, on the other hand, $\mathcal{M}, w,h \vDash \Box B \to B $, then $\Box B \to B$ is true at all points of $U$. Consequently, $\Box (\Box B \to B)$ is true at all points of $U$.
Applying (\ref{form1}), we obtain
\begin{gather*}
\forall w^\prime\in U\; (\rho_\tau(w^\prime) < \rho_\tau(w) \Longrightarrow \mathcal{M}, w^\prime, h \vDash  \Box B).
\end{gather*}
Recalling that $\Box B \to B$ is true at all points of $U$, we also have
\[\forall w^\prime\in U\; (\rho_\tau(w^\prime) < \rho_\tau(w) \Longrightarrow \mathcal{M}, w^\prime, h \vDash B).\]
Since $\{w^\prime\in X \mid \rho_\tau(w^\prime) < \rho_\tau(w)\}$ is a punctured neighborhood of $w$ by Lemma \ref{punctured neighbourhood} and $U$ is a neighborhood of $w$, the set $ \{w^\prime\in U \mid\rho_\tau(w^\prime) < \rho_\tau(w) \}= U\cap \{w^\prime\in X \mid \rho_\tau(w^\prime) < \rho_\tau(w)\}$ is a punctured neighborhood of $w$. Moreover, the formula $B$ is true at every point of this set. We obtain that $\mathcal{M}, w,h \vDash \Box B  $ and $\mathcal{M}, w,h \vDash A $.

Case 2. Suppose that the $\infty$-proof $\pi$ does not consist only of the root, and consider the last application of an inference rule in $\pi$.

Case 2.1. Let $\pi$ have the form
\begin{gather*}
\AXC{$\pi^\prime$}
\noLine
\UIC{$\vdots$}
\noLine
\UIC{$B$}
\AXC{$\pi^{\prime\prime}$}
\noLine
\UIC{$\vdots$}
\noLine
\UIC{$B \to A$}
\LeftLabel{$\mathsf{mp}$}
\RightLabel{ .}
\BIC{$A$}
\DisplayProof
\end{gather*}
Since $\lvert  \pi^\prime\rvert< \lvert  \pi\rvert$ and $\lvert  \pi^{\prime\prime}\rvert< \lvert  \pi\rvert$, by the induction hypothesis $\mathcal{M}, w, h \vDash B$ and $\mathcal{M}, w, h \vDash B \to A$. Consequently, $\mathcal{M}, w, h \vDash A$.

Case 2.2. Now suppose that $\pi$ has the form
\begin{gather*}
\AXC{$\pi^\prime$}
\noLine
\UIC{$\vdots$}
\noLine
\UIC{$B$}
\LeftLabel{$\mathsf{gen}$}
\RightLabel{ ,}
\UIC{$\forall x\: B$}
\DisplayProof
\end{gather*}
where $\forall x\: B =A$. Note that $\lvert  \pi^\prime\rvert< \lvert  \pi\rvert$. By the induction hypothesis applied to $h^a_x$, where $a$ is an arbitrary element of $D$, we have $\mathcal{M}, w, h^a_x \vDash B$. Consequently, $\mathcal{M}, w, h \vDash \forall x\: B$.

Case 2.3. It remains to consider the case when $\pi$ has the form
\begin{gather*}
\AXC{$\pi^\prime$}
\noLine
\UIC{$\vdots$}
\noLine
\UIC{$B$}
\LeftLabel{$\mathsf{nec}$}
\RightLabel{ ,}
\UIC{$\Box B$}
\DisplayProof
\end{gather*}
where $\Box B =A$. By the induction hypothesis applied to $\pi^\prime$,
\[\forall w^\prime\in X\; (\rho_\tau(w^\prime) < \rho_\tau(w) \Longrightarrow \mathcal{M}, w^\prime, h \vDash B).\]
Since $\{w^\prime\in X \mid\rho_\tau(w^\prime) < \rho_\tau(w)\}$ is a punctured neighborhood of $w$ at all points of which the formula $B$ is true, we obtain $\mathcal{M}, w, h \vDash \Box B$.




\end{proof}

\section{Sequent calculus}
\label{s3}
In this section, we define a sequent calculus $\mathsf{S}$ admitting non-well-founded proofs, and establish completeness of the system $\mathsf{S}$ with respect to topological semantics. Relying on the latter result, in the next section we prove the topological completeness of the system $\mathsf{QGL}_\infty$.


Recall that a \emph{sequent} is any expression of the form $\G \Rightarrow \D$, where $\Gamma$ and $\Delta$ are finite multisets of formulas. 
The \emph{formula image} of a sequent $\G \Rightarrow \D$ is the formula ${\bigwedge\G \to \bigvee\D}$. Here, $\bigwedge \emptyset \coloneq \top$ and $\bigvee \emptyset \coloneq \bot$.


For a finite multiset of formulas $\Gamma = \{A_0, \dotsc , A_{n-1}\}$, we define the multiset $\Box \Gamma \coloneq \{\Box A_0, \dotsc , \Box A_{n-1}\}$. Also let $\Gamma^s$ be the set of elements of $\Gamma$. For a finite set of formulas $S$, we define $S^m$ as the multiset that contains exactly all elements of $S$ and, moreover, without repetitions.

The union of multisets is denoted by the symbol $\cupplus$. However, we will often write $\Gamma, \Delta$ instead of $\Gamma \cupplus \Delta$, and also omit the curly braces when writing multisets.

The initial sequents and inference rules of the calculus $\mathsf{S}$ have the following form:
\begin{gather*}
\G, P(y_0,\dotsc, y_{n-1})\Rto P(y_0,\dotsc, y_{n-1}), \D \: , \qquad \G ,\bot \Rto \Delta \:,\\\\
\AXC{$\G, B \Rto \D$}
\AXC{$\G\Rto A, \D$}
\LeftLabel{$\mathsf{\to_L}$}
\RightLabel{,}
\BIC{$\G, A\to B \Rto \D$}
\DisplayProof
\qquad
\AXC{$\G, A \Rto B, \D$}
\LeftLabel{$\mathsf{\to_R}$}
\RightLabel{,}
\UIC{$\G \Rto A \to B, \D$}
\DisplayProof\\\\
\AXC{$\G, A[x\mapsto y], \forall x \: A \Rto \D$}
\LeftLabel{$\mathsf{\forall_L}$}
\RightLabel{,}
\UIC{$\G, \forall x \: A \Rto \D$}
\DisplayProof
\qquad
\AXC{$\G \Rto \D, A[x\mapsto y]$}
\LeftLabel{$\mathsf{\forall_R}$}
\RightLabel{($y \notin \mathit{FV}(\G , \D ,\forall x \: A)$),}
\UIC{$\G \Rto \D, \forall x \: A$}
\DisplayProof\\\\
\AXC{$\G, \Box\G \Rto A$}
\LeftLabel{$\Box$}
\RightLabel{.}
\UIC{$\Pi ,\Box\G \Rto \Box A, \D$}
\DisplayProof
\end{gather*}
Note that for a multiset of formulas $\Pi$, the notation $\mathit{FV}(\Pi)$ means the set of all individual variables having free occurrences in formulas from $\Pi$.

A \emph{non-well-founded proof, or $\infty$-proof, in the sequent calculus $\mathsf{S}$} is a (possibly infinite) tree whose vertices are labeled with sequents, whose leaves are labeled with initial sequents, and where the transition from children to parents is performed according to the inference rules of the calculus $\mathsf{S}$. Moreover, every infinite branch of an $\infty$-proof must intersect the rule ($\Box$) infinitely many times. A sequent $\Gamma\Rightarrow \Delta$ is \emph{provable in $\mathsf{S}$} if there exists an $\infty$-proof whose root is labeled with the sequent $\Gamma\Rightarrow \Delta$.


Similarly to the case of $\infty$-proofs in the calculus $\mathsf{QGL}_\infty$, we define the \emph{main fragment of an $\infty$-proof $\pi$ in the calculus $\mathsf{S}$} as the tree obtained from $\pi$ by cutting each branch at the application of the rule ($\Box$) closest to the root. The \emph{local height $\lvert \pi \rvert$ of an $\infty$-proof $\pi$} is the length of the longest branch of the main fragment of $\pi$. If $\pi$ consists of a single vertex, then $\lvert \pi \rvert=0$.

We call a pair $(L,R)$, where $L$ and $R$ are sets of formulas, \emph{$\mathit{CB}$-valid} if for every $\mathsf{QGL}$-model $\mathcal{M}=(X,\tau, D, \xi)$, every world $w\in X$, and every assignment $h\colon \mathit{Var} \to D$,
\[\forall B\in L \: \mathcal{M}, w, h \vDash B \Longrightarrow \exists C\in R \: \mathcal{M}, w, h \vDash C,\]
provided that the domain $D$ is countable and, moreover, $h\colon \mathit{Var} \to D$ is bijective. A sequent $\Gamma\Rightarrow \Delta$ is \emph{$\mathit{CB}$-valid} if the pair $(\Gamma^s, \Delta^s)$ is $\mathit{CB}$-valid.


\begin{lem}\label{semantic inversion}
If the conclusion of one of the inference rules $(\mathsf{\to_L})$, $(\mathsf{\to_R})$, $(\mathsf{\forall_L})$, or $(\mathsf{\forall_R})$ is $\mathit{CB}$-valid, then the premises of the rule are also $\mathit{CB}$-valid.
\end{lem}

A pair $(L,R)$, where $L$ and $R$ are sets of formulas, is called \emph{saturated} if it satisfies the following conditions:
\begin{itemize}
\item $\bot\notin L$;
\item $A\in R$ or $B\in L$, if $(A\to B) \in L$;
\item $A\in L$ and $B\in R$, if $(A\to B) \in R$;
\item $A[x\mapsto y]\in L$ for every $y\in \mathit{Var}$, if $\forall x\: A \in L$;
\item $A[x\mapsto y]\in R$ for some $y\in \mathit{Var}$, if $\forall x\: A \in R$.
\end{itemize}

\begin{lem}\label{saturation lemma}
Let a pair $(L,R)$, where $L$ and $R$ are sets of formulas, be saturated and $\mathit{CB}$-valid. Then either there is an atomic formula $P(y_0,\dotsc, y_{n-1}) \in L \cap R$, or there are a finite set of formulas $S$ and a formula $A$ such that $\Box S\subset L$, $\Box A \in R$, and the sequent $S^m, \Box S^m \Rightarrow A$ is $\mathit{CB}$-valid.
\end{lem}
\begin{proof}
Let $(L,R)$ be a saturated and $\mathit{CB}$-valid pair of sets of formulas. We establish the statement of the lemma by contradiction.

Suppose that there is no atomic formula that belongs simultaneously to both $L$ and $R$. Suppose also that for every finite set of formulas $S$ and every formula $A$ such that $\Box S \subset L$ and $\Box A\in R$, the sequent $S^m, \Box S^m\Rightarrow A$ is not $\mathit{CB}$-valid.

Let $I\coloneq \{(S,A)\in \mathcal{P}(\mathit{Fm}_\Omega)\times  \mathit{Fm}_\Omega\mid \text{$S$ is finite, $\Box S \subset L$ and $\Box A\in R$}\}$. Furthermore, let $\mathit{pr}_0(i)$ and $\mathit{pr}_1(i)$ denote the first and second components of a pair $i\in I$.
By the assumptions made above, for each $i\in I$ we can fix a $\mathsf{QGL}$-model $\mathcal{M}_i=(X_i,\tau_i, D_i, \xi_i)$, a world $w_i\in X_i$, and a bijective assignment $h_i\colon \mathit{Var} \to D_i$ such that
\[\mathcal{M}_i, w_i, h_i \vDash \bigwedge S\wedge\bigwedge\Box S \qquad \text{ and } \qquad\mathcal{M}_i, w_i, h_i \nvDash A,\]
where $S= \mathit{pr}_0(i)$ and $A= \mathit{pr}_1(i)$.
Since for each $i\in I$ the assignment $h_i$ is bijective, without loss of generality we may assume that all domains $D_i$ are equal to each other and coincide with the set $D=\{a_0, a_1, a_2, \dotsc\}$, and the assignments $h_i$ coincide with the function $h\colon x_i \mapsto a_i$. Also, for $i,j\in I$, we assume without loss of generality that $X_i\cap X_j=\emptyset$ if $i\neq j$.


Let us define a $\mathsf{QGL}$-frame that will allow us to refute the $\mathit{CB}$-validity of the pair $(L,R)$. Put $X \coloneq \{r\}\cup \bigcup\limits_{i\in I} X_i$, where $r$ is a new element not belonging to $\bigcup\limits_{i\in I} X_i$, and introduce a topology $\tau$ on $X$. Let a set $U\in \tau$ if and only if $U\subset X$ and
\begin{itemize}
\item $U\cap X_i \in \tau_i$ for all $i\in I$,
\item if $r\in U$, then there is a finite set of formulas $S$ such that $\Box S \subset L$ and, moreover, $w_i\in U$ whenever $S\subset \mathit{pr}_0(i)$.
\end{itemize}
It is obvious that $(X, \tau)$ is a scattered topological space and $(X, \tau, D)$ is a $\mathsf{QGL}$-frame.

Now we define an indexed family $\alpha=(\alpha_w)_{w\in X}$ of valuations in $D$. If the signature contains a predicate symbol $P$ of arity $n$, then
\begin{itemize}
\item $\alpha_w (P)\coloneq(\xi_i)_w (P)$ for $w\in X_i$,
\item $\alpha_r (P)\coloneq\{(h(y_0),\dotsc, h(y_{n-1})) \in D^n\mid P(y_0,\dotsc, y_{n-1})\in L\}$.
\end{itemize}

We have constructed a $\mathsf{QGL}$-model $\mathcal{M}=(X, \tau, D, \alpha)$, with the help of which it is easy to refute the $\mathit{CB}$-validity of the pair $(L,R)$. Indeed, let us check that for every formula $B$ the following holds:
\begin{gather}\label{form2}
B\in L\Longrightarrow \mathcal{M}, r, h \vDash B \qquad \text{ and } \qquad B\in R\Longrightarrow\mathcal{M}, r, h \nvDash B.
\end{gather}
We establish (\ref{form2}) by induction on the length of the formula $B$.

Case 1. Suppose that $B=\bot$. Note that $\bot\notin L$ by the saturation of the pair $(L,R)$. At the same time, $\mathcal{M}, r, h \nvDash \bot$ by the definition of the truth relation.

Case 2. Let $B$ be of the form $P(y_0,\dotsc, y_{n-1})$. Then by the definition of the valuation $\alpha_r$,
\[P(y_0,\dotsc, y_{n-1})\in L\Longleftrightarrow \mathcal{M}, r, h \vDash P(y_0,\dotsc, y_{n-1}).\]
Moreover,
\[P(y_0,\dotsc, y_{n-1})\in R \Longrightarrow P(y_0,\dotsc, y_{n-1})\notin L\Longrightarrow \mathcal{M}, r, h \nvDash P(y_0,\dotsc, y_{n-1}),\]
since there is no atomic formula that belongs simultaneously to both $L$ and $R$.

Case 3. Suppose that $B=B'\to B''$. By the saturation of the pair $(L,R)$ and by the induction hypothesis,
\begin{multline*}
B\in L\Longrightarrow (B'\in R \text{ or } B''\in L)\Longrightarrow\\ \Longrightarrow (\mathcal{M}, r, h \nvDash B' \text{ or } \mathcal{M}, r, h \vDash B'') \Longrightarrow \mathcal{M}, r, h \vDash B.
\end{multline*}
Also,
\begin{multline*}
B\in R\Longrightarrow (B'\in L \text{ and } B''\in R)\Longrightarrow\\ \Longrightarrow (\mathcal{M}, r, h \vDash B' \text{ and } \mathcal{M}, r, h \nvDash B'') \Longrightarrow \mathcal{M}, r, h \nvDash B.
\end{multline*}

Case 4. Let $B$ be of the form $\Box B'$. Suppose that $B\in L$. Define $J\coloneq\{i\in I\mid B'\in \mathit{pr}_0(i)\}$ and note that $\mathcal{M}_j, w_j, h \vDash B'\wedge \Box B'$ for each $j\in J$. By the definition of truth of a formula in a $\mathsf{QGL}$-model, for every $j\in J$ there is $U_j\in \tau_j$ such that $w_j \in U_j$ and $\mathcal{M}_j, w, h \vDash B'$ for each $w\in U_j$. It follows by Lemma \ref{submodel} that $\mathcal{M}, w, h \vDash B'$ for every $j\in J$ and every $w\in U_j$. Thus, $\mathcal{M}, w, h \vDash B'$ for all $w\in \bigcup\limits_{j\in J} U_j$. Moreover, by the definition of the topology $\tau$ on $X$, the set $\bigcup\limits_{j\in J} U_j$ is a punctured neighborhood of $r$. Hence, $\mathcal{M}, r, h \vDash \Box B'$.

Now suppose that $B\in R$, and consider an arbitrary open set $U$ containing $r$. By the definition of $\tau$, there is a finite set of formulas $S$ such that $\Box S \subset L$ and, moreover, $w_i\in U$ whenever $S\subset \mathit{pr}_0(i)$. Note that $(S, B')\in I$ and $\mathcal{M}_{(S, B')}, w_{(S, B')}, h \nvDash B'$. By Lemma \ref{submodel}, $\mathcal{M}, w_{(S, B')}, h \nvDash B'$. Here, $w_{(S, B')}\in U\setminus \{r\}$. Thus, $\mathcal{M}, r, h \nvDash \Box B'$.

Case 5. It remains to consider the case when $B$ is of the form $\forall x\: B'$. By the saturation of the pair $(L,R)$ and by the induction hypothesis,
\begin{multline*}
B\in L\Longrightarrow \forall y\in \mathit{Var} \:B'[x\mapsto y]\in L \Longrightarrow \forall y\in \mathit{Var} \:\mathcal{M}, r, h \vDash B'[x\mapsto y]\Longrightarrow \\ \Longrightarrow \forall a\in D\: \mathcal{M}, r, h^a_x \vDash B' \Longrightarrow \mathcal{M}, r, h \vDash B.
\end{multline*}
At the same time,
\begin{multline*}
B\in R\Longrightarrow \exists y\in \mathit{Var} \:B'[x\mapsto y]\in R \Longrightarrow \exists y\in \mathit{Var} \:\mathcal{M}, r, h \nvDash B'[x\mapsto y]\Longrightarrow \\ \Longrightarrow \exists a\in D\: \mathcal{M}, r, h^a_x \nvDash B' \Longrightarrow \mathcal{M}, r, h \nvDash B.
\end{multline*}

Thus, condition (\ref{form2}) is established. From this it immediately follows that the pair $(L,R)$ is not $\mathit{CB}$-valid. We have obtained a contradiction with the initial assumption. Therefore, either there is an atomic formula $P(y_0,\dotsc, y_{n-1}) \in L \cap R$, or there are a finite set of formulas $S$ and a formula $A$ such that $\Box S\subset L$, $\Box A \in R$, and the sequent $S^m, \Box S^m \Rightarrow A$ is $\mathit{CB}$-valid.

\end{proof}
The following completeness statement for the system $\mathsf{S}$ will play a key role in the proof of the topological completeness theorem for the calculus $\qgl_{\infty}$.
\begin{prop}\label{cb-completeness}
If a sequent $\G \Rto \D$ is $\mathit{CB}$-valid, then $\mathsf{S} \vdash \G \Rto \D$.
\end{prop}
\begin{proof}
Consider a $\mathit{CB}$-valid sequent $\G \Rto \D$ and find an $\infty$-proof for it in the calculus $\mathsf{S}$ using a variation of K. Sch\"utte's method based on the construction of reduction trees (see \cite{schutte1956, takeuti1978, negri2001}). In other words, the desired $\infty$-proof is obtained by applying the following algorithm to the tree $\delta$ consisting of a single vertex labeled with the sequent $\G \Rto \D$.





Action 1. If the current tree $\delta$ has a leaf labeled with a sequent that is not initial and that can be obtained by applying the rule ($\mathsf{\to_L}$), then we extend $\delta$ by placing the two premises of the rule ($\mathsf{\to_L}$) above this leaf and repeat Action 1. If there is no such leaf, we proceed to Action 2.

Action 2. If the current tree $\delta$ has a leaf labeled with a sequent that is not initial and that can be obtained by applying the rule ($\mathsf{\to_R}$), then we extend $\delta$ by placing the premise of the rule ($\mathsf{\to_R}$) above this leaf and repeat Action 2. If there is no such leaf, we proceed to Action 3.

Action 3. Above each leaf $e$ of the tree $\delta$ labeled with a sequent $\Pi, \forall y_0\: A_0, \dotsc ,\forall y_{n}\: A_{n} \Rightarrow \Lambda $ that is not initial and in which $\Pi$ contains no formulas of the form $\forall y\: A$, we place the premise of the rule ($\mathsf{\forall_L}$) $n+1$ times so that the new leaf is labeled with the sequent $\Pi, A_0 [y_0\mapsto x_{i_0}], \forall y_0\: A_0, \cdots , A_n [y_n\mapsto x_{i_n}], \forall y_n\: A_n  \Rightarrow \Lambda$. Here, $i_k$ for each $k\in \{0, \dotsc , n \}$ is defined as the least natural number for which the formula $A_k [y_k \mapsto x_{i_k}]$ does not occur in the left-hand sides of sequents lying between the leaf $e$ and the nearest premise of the rule ($\Box$) below, including the leaf and the premise, or between this leaf $e$ and the root if there is no such premise. After this, we proceed to Action 4.

Action 4. If the current tree $\delta$ has a leaf labeled with a sequent $\Pi \Rightarrow \Lambda, \forall x\: A $ that is not initial, then we place above this leaf the sequent $\Pi \Rightarrow \Lambda, A [x\mapsto x_i]$, where $i=\min \{j\in \mathbb{N}\mid x_j \notin \mathit{FV}(\Pi , \Lambda ,\forall x \: A))\}$, and repeat Action 4. If there is no such leaf, we proceed to Action 5.

Action 5. Above each leaf of the tree $\delta$ labeled with a sequent that is not initial and that can be obtained by applying the rule ($\Box$) from a $\mathit{CB}$-valid premise, we place such a premise. After this, we proceed to Action 1.

The described algorithm, in the process of its application, constructs a tree $\pi$ whose vertices are labeled with sequents, the transition from children to parents is performed according to the inference rules of the calculus $\mathsf{S}$, and the root is labeled with the sequent $\Gamma \Rightarrow \Delta$. By the definition of the algorithm and by Lemma \ref{semantic inversion}, the obtained tree $\pi$ consists entirely of $\mathit{CB}$-valid sequents.

Now let us verify that $\pi$ is an $\infty$-proof in the calculus $\mathsf{S}$. To do this, it suffices to check that every leaf of $\pi$ is labeled with an initial sequent, and that every infinite branch intersects the rule ($\Box$) infinitely many times.

Consider an arbitrary leaf $e$ of the tree $\pi$ and the sequent $\Pi \Rightarrow \Lambda$ with which the leaf $e$ is labeled. Arguing by contradiction, suppose that $\Pi \Rightarrow \Lambda$ is not initial. Then $\Pi$ does not contain the formula $\bot$. Moreover, the multiset $\Pi$ (the multiset $\Pi \cupplus \Lambda$) contains no formulas of the form $\forall y\: A$ (of the form $A\to B$); otherwise Action 3 (Action 1 or 2) would be applicable to $e$. Consequently, the pair $(\Pi^s, \Lambda^s)$ is saturated and, moreover, $\mathit{CB}$-valid. It follows from Lemma \ref{saturation lemma} that $\Pi \Rightarrow \Lambda$ can be obtained from a $\mathit{CB}$-valid premise by means of the rule ($\Box$). Thus, Action 5 is applicable to $e$, and $e$ cannot be a leaf of the tree $\pi$. This contradiction shows that $\Pi \Rightarrow \Lambda$ is an initial sequent.

Now consider an arbitrary infinite branch of the tree $\pi$
\[e_0 \to e_1 \to e_2 \to  \dotsb\]
and check that it intersects the rule ($\Box$) infinitely many times, arguing by contradiction. Suppose that this branch intersects the rule ($\Box$) only finitely many times. Then it has a tail
\[e_m \to e_{m+1} \to e_{m+2} \to  \dotsb,\]
which does not intersect the rule ($\Box$), where $e_m$ is either the root of the tree $\pi$ or a premise of the rule ($\Box$). Denote by $\Pi_i\Rightarrow \Lambda_i$ the sequents with which the vertices $e_{m+i}$ are labeled. By the construction of the tree $\pi$, the sequence $(\Pi_i\Rightarrow \Lambda_i)_{i\in \mathbb{N}}$ contains no initial sequents. In particular, $\bot\notin \Pi_i$ for each $i\in \mathbb{N}$.
Then by the construction of $\pi$, the pair $(\bigcup\limits_{i\in \mathbb{N}}  \Pi^s_i, \bigcup\limits_{i\in \mathbb{N}}\Lambda^s_i)$ is saturated and, moreover, $\mathit{CB}$-valid, since the sequent $\Pi_0\Rightarrow \Lambda_0$ is $\mathit{CB}$-valid.
From this, by Lemma \ref{saturation lemma}, we obtain that either there is an atomic formula $P(y_0,\dotsc, y_{n-1}) \in \bigcup\limits_{i\in \mathbb{N}}  \Pi^s_i \cap \bigcup\limits_{i\in \mathbb{N}}\Lambda^s_i$, or there are a finite set of formulas $S$ and a formula $A$ such that $\Box S\subset \bigcup\limits_{i\in \mathbb{N}}  \Pi^s_i$, $\Box A \in \bigcup\limits_{i\in \mathbb{N}}\Lambda^s_i$, and the sequent $S^m, \Box S^m \Rightarrow A$ is $\mathit{CB}$-valid. In the first case,
among the sequents $\Pi_i\Rightarrow \Lambda_i$ there is an initial one, which is impossible. In the second case, Action 5 is applicable to all vertices of the tail starting from some point. Then by the construction of the tree $\pi$, the tail intersects the rule ($\Box$), which contradicts the choice of the tail. Thus, the original branch intersects the rule ($\Box$) infinitely many times.

\end{proof}

\section{Topological completeness}
\label{s4}

\begin{prop}\label{proof transformation}
If $\mathsf{S} \vdash \G \Rto \D$, then $\mathsf{QGL}_\infty \vdash \bigwedge\G \to \bigvee\D$.
\end{prop}
\begin{proof}
Consider an $\infty$-proof $\pi$ of the sequent $\Gamma \Rightarrow \Delta$ in the calculus $\mathsf{S}$. We define a transformation of this $\infty$-proof into the desired $\infty$-proof $f(\pi)$ of the system $\mathsf{QGL}_\infty$, moving upward through the tree $\pi$ and passing in the inference rules from conclusions to premises. More precisely, we define the transformation $f$ by co-recursion with a nested recursion on the local height $\lvert\pi \rvert$. For a detailed discussion of such definitions, we refer the reader to \cite{miranda2024}.

Case 1. If $\pi$ consists of a single vertex, then $\Gamma \Rightarrow \Delta$ is an initial sequent. In this case, the formula $\bigwedge \Gamma \rightarrow \bigvee \Delta$ is provable in $\mathsf{QGL}_\infty$ by a finite $\infty$-proof. Let $f(\pi)$ be such a proof.

Case 2. Suppose that the $\infty$-proof $\pi$ consists of more than one vertex, and consider the last application of an inference rule in $\pi$.

Case 2.1. Suppose that $\pi$ has the form
\begin{gather*}
\AXC{$\pi^\prime$}
\noLine
\UIC{$\vdots$}
\noLine
\UIC{$ \Pi, B\Rightarrow \Delta$}
\AXC{$\pi^{\prime\prime}$}
\noLine
\UIC{$\vdots$}
\noLine
\UIC{$ \Pi \Rightarrow A,\Delta$}
\LeftLabel{$\mathsf{\to_L}$}
\RightLabel{ ,}
\BIC{$ \Pi, A\to B\Rightarrow \Delta$}
\DisplayProof
\end{gather*}
where $ \Pi \cupplus\{A\to B\}=\Gamma$. Since $\lvert  \pi^\prime\rvert< \lvert  \pi\rvert$ and $\lvert  \pi^{\prime\prime}\rvert< \lvert  \pi\rvert$, the $\infty$-proofs $f(\pi^\prime)$ and $f(\pi^{\prime\prime})$ are defined. Let $f(\pi)$ be equal to
\begin{gather*}
\AXC{$f(\pi^\prime)$}
\noLine
\UIC{$\vdots$}
\noLine
\UIC{$G$}
\AXC{$f(\pi^{\prime\prime})$}
\noLine
\UIC{$\vdots$}
\noLine
\UIC{$F$}
\AXC{$\eta$}
\noLine
\UIC{$\vdots$}
\noLine
\UIC{$F \rightarrow (G \rightarrow H)$}
\LeftLabel{$\mathsf{mp}$}
\BIC{$G\to H$}
\LeftLabel{$\mathsf{mp}$}
\RightLabel{ ,}
\BIC{$H$}
\DisplayProof
\end{gather*}
where $F=\bigwedge\Pi \rightarrow \bigvee(\{A\} \cupplus \Delta)$, $G=\bigwedge (\Pi \cupplus \{B\} )\rightarrow \bigvee\Delta$, $H= \bigwedge (\Pi\cupplus \{A\to B\})\rightarrow \bigvee\Delta$, and $\eta$ is a finite $\infty$-proof of the formula $F \rightarrow (G \rightarrow H)$ in $\mathsf{QGL}_\infty$.

Case 2.2. Suppose that $\pi$ has the form
\begin{gather*}
\AXC{$\pi^\prime$}
\noLine
\UIC{$\vdots$}
\noLine
\UIC{$\Gamma,A\Rightarrow B,\Lambda$}
\LeftLabel{$\mathsf{\to_R}$}
\RightLabel{ ,}
\UIC{$\Gamma\Rightarrow A\to B, \Lambda$}
\DisplayProof
\end{gather*}
where $\{A\to B\}\cupplus \Lambda =\Delta$. Since $\lvert  \pi^\prime\rvert< \lvert  \pi\rvert$, the $\infty$-proof $f(\pi^\prime)$ of the formula $\bigwedge (\Gamma\cupplus \{A\}) \rightarrow \bigvee(\{B\} \cupplus \Lambda)$ is defined. Let $f(\pi)$ be equal to
\begin{gather*}
\AXC{$f(\pi^{\prime})$}
\noLine
\UIC{$\vdots$}
\noLine
\UIC{$F$}
\AXC{$\eta$}
\noLine
\UIC{$\vdots$}
\noLine
\UIC{$F \rightarrow G$}
\LeftLabel{$\mathsf{mp}$}
\RightLabel{ ,}
\BIC{$G$}
\DisplayProof
\end{gather*}
where $F=\bigwedge (\Gamma\cupplus \{A\}) \rightarrow \bigvee(\{B\} \cupplus \Lambda)$, $G= \bigwedge \Gamma\rightarrow \bigvee(\{A\to B\}\cupplus \Lambda)$, and $\eta$ is a finite $\infty$-proof of the formula $F \rightarrow G$ in $\mathsf{QGL}_\infty$.

Case 2.3. Consider the case when $\pi$ has the form
\begin{gather*}
\AXC{$\pi^\prime$}
\noLine
\UIC{$\vdots$}
\noLine
\UIC{$\Pi, A[x\mapsto y], \forall x \: A \Rto \D$}
\LeftLabel{$\mathsf{\forall_L}$}
\RightLabel{,}
\UIC{$\Pi, \forall x \: A \Rto \D$}
\DisplayProof
\end{gather*}
where $ \Pi \cupplus\{\forall x \: A\}=\Gamma$. Since $\lvert  \pi^\prime\rvert< \lvert  \pi\rvert$, the $\infty$-proof $f(\pi^\prime)$ of the formula $\bigwedge (\Pi\cupplus \{A[x\mapsto y], \forall x \: A \}) \rightarrow \bigvee \Delta$ is defined. Define $f(\pi)$ as
\begin{gather*}
\AXC{$f(\pi^{\prime})$}
\noLine
\UIC{$\vdots$}
\noLine
\UIC{$F$}
\AXC{$\eta$}
\noLine
\UIC{$\vdots$}
\noLine
\UIC{$F \rightarrow G$}
\LeftLabel{$\mathsf{mp}$}
\RightLabel{ ,}
\BIC{$G$}
\DisplayProof
\end{gather*}
where $F=\bigwedge (\Pi\cupplus \{A[x\mapsto y], \forall x \: A \}) \rightarrow \bigvee \Delta$, $G= \bigwedge (\Pi\cupplus \{ \forall x \: A \}) \rightarrow \bigvee \Delta$, and $\eta$ is a finite $\infty$-proof of the formula $F \rightarrow G$ in $\mathsf{QGL}_\infty$.

Case 2.4. Let $\pi$ have the form
\begin{gather*}
\AXC{$\pi^\prime$}
\noLine
\UIC{$\vdots$}
\noLine
\UIC{$\G \Rto  A[x\mapsto y], \Lambda$}
\LeftLabel{$\mathsf{\forall_R}$}
\RightLabel{($y \notin FV(\G , \Lambda , \forall x \: A)$),}
\UIC{$\G \Rto  \forall x \: A, \Lambda$}
\DisplayProof
\end{gather*}
where $\{ \forall x \: A \} \cupplus \Lambda = \Delta$. The $\infty$-proof $f(\pi^\prime)$ of the formula $\bigwedge \Gamma \rightarrow \bigvee (\{A[x\mapsto y]\} \cupplus \Lambda)$ is defined, because $\lvert  \pi^\prime\rvert< \lvert  \pi\rvert$. Let $f(\pi)$ be equal to
\begin{gather*}
\AXC{$f(\pi^{\prime})$}
\noLine
\UIC{$\vdots$}
\noLine
\UIC{$F$}
\LeftLabel{$\mathsf{gen}$}
\UIC{$\forall y \: F$}
\AXC{$\eta$}
\noLine
\UIC{$\vdots$}
\noLine
\UIC{$\forall y \: F \rightarrow G$}
\LeftLabel{$\mathsf{mp}$}
\RightLabel{ ,}
\BIC{$G$}
\DisplayProof
\end{gather*}
where $F=\bigwedge \Gamma \rightarrow \bigvee (\{A[x\mapsto y]\} \cupplus \Lambda)$, $G= \bigwedge \Gamma \rightarrow \bigvee (\{\forall x \: A\} \cupplus \Lambda)$, and $\eta$ is a finite $\infty$-proof of the formula $\forall y \: F \rightarrow G$ in $\mathsf{QGL}_\infty$.

Case 2.5. It remains to consider the case when $\pi$ has the form
\begin{gather*}
\AXC{$\pi^\prime$}
\noLine
\UIC{$\vdots$}
\noLine
\UIC{$ \Upsilon, \Box \Upsilon \Rightarrow A$}
\LeftLabel{$\mathsf{\Box}$}
\RightLabel{,}
\UIC{$\Pi ,\Box \Upsilon \Rightarrow \Box A , \Lambda$}
\DisplayProof
\end{gather*}
where $ \Pi \cupplus \Box \Upsilon=\Gamma$ and $ \{\Box A\} \cupplus \Lambda=\Delta$. Note that $\lvert  \pi^\prime\rvert=0$. Using co-recursion, we apply the transformation $f$ to the $\infty$-proof $\pi^\prime$ and define $f(\pi)$ as
\begin{gather*}
\AXC{$f(\pi^\prime)$}
\noLine
\UIC{$\vdots$}
\noLine
\UIC{$F$}
\LeftLabel{$\mathsf{nec}$}
\UIC{$\Box F$}
\AXC{$\eta$}
\noLine
\UIC{$\vdots$}
\noLine
\UIC{$\Box F \rightarrow G$}
\LeftLabel{$\mathsf{mp}$}
\RightLabel{ ,}
\BIC{$G$}
\DisplayProof
\end{gather*}
where $F=\bigwedge (\Upsilon\cupplus \Box \Upsilon) \rightarrow A$, $G= \bigwedge (\Pi\cupplus \Box \Upsilon) \rightarrow \bigvee(\{\Box A\} \cupplus \Lambda)$, and $\eta$ is a finite $\infty$-proof of the formula $\Box F \rightarrow G$ in $\mathsf{QGL}_\infty$.

The definition of the transformation $f$ is complete. The obtained $f(\pi)$ is a tree whose vertices are labeled with $\mathsf{QGL}$-formulas, and the transition from children to parents is performed according to the inference rules ($\mathsf{mp}$), ($\mathsf{gen}$), and ($\mathsf{nec}$). Here, the leaves of this tree are labeled with axioms, and every infinite branch intersects the rule ($\mathsf{nec}$) infinitely many times. Thus, $f(\pi)$ is the desired $\infty$-proof. We have established that $\mathsf{QGL}_\infty \vdash \bigwedge\G \to \bigvee\D$.
\end{proof}

\begin{thm}[topological completeness]\label{completeness}
A formula $A$ is valid if and only if $\mathsf{QGL}_\infty \vdash A$.
\end{thm}

\begin{proof}
It suffices to show that if a formula $A$ is valid, then $\mathsf{QGL}_{\infty} \vdash A$, since the converse holds by Proposition \ref{soundness}. If the formula $A$ is valid, then the sequent $\emptyset \Rto A$ is $\mathit{CB}$-valid. By Proposition \ref{cb-completeness}, the sequent $\emptyset \Rto A$ is provable in $\mathsf{S}$. Then $\mathsf{QGL}_{\infty} \vdash \top \to A$ by Proposition \ref{proof transformation}. From this we immediately obtain $\mathsf{QGL}_{\infty} \vdash A$.
\end{proof}

\section{Conclusion}
\label{s5}
The established result on the topological completeness of the calculus $\mathsf{QGL}_\infty$ leads to new questions both in the area of modal predicate logics and systems with non-well-founded derivations, and in the area of predicate provability logic. Let us note several directions for possible research.

First, the question of strong global completeness of the calculus $\mathsf{QGL}_\infty$ with respect to topological semantics is of interest. It is known that in the case of the logic $\mathsf{GL}$ extended by non-well-founded derivations, global completeness holds \cite{shamkanov2020}. At the same time, whether this result holds in the predicate case remains unclear. It is possible that a proof of global completeness in the case of $\mathsf{QGL}_\infty$ will require a transition from topological semantics to a more general one, to the semantics of (pre)sheaves.

Second, there arises the problem of obtaining an exact complexity bound for the set of theorems of the calculus $\mathsf{QGL}_\infty$, and a similar problem for the variant of the calculus $\mathsf{QGL}_\infty$ in which only computable non-well-founded derivations are allowed. Recently, one of the above-mentioned problems has apparently been solved.

Finally, the question of the arithmetical semantics of the calculus $\mathsf{QGL}_\infty$ is of considerable interest. The deep investigation \cite{pakhomov2021} opens an approach to proving arithmetical soundness in the case when, instead of the arithmetic $\mathsf{PA}$, we take various second-order systems, for example, the arithmetic $\mathsf{ACA}_0$ extended by all true $\Sigma^1_1$-sentences. At the same time, the prospects for proving arithmetical completeness are not yet clear.

\end{document}